\documentclass[a4paper,12pt]{article}
\usepackage[T2A]{fontenc}
\usepackage[cp1251]{inputenc}
\usepackage{amsmath}
\usepackage{amssymb}
\usepackage[dvips]{graphicx}
\usepackage[english,russian]{babel}

\newtheorem{theorem}{Теорема}[section]
\newtheorem{corollary}[theorem]{Следствие}

\newtheorem{proposition}[theorem]{Предложение}
\newtheorem{definitionhead}[theorem]{Определение}
\newenvironment{definition}{\begin{definitionhead}%
\sl}{\end{definitionhead}}

\newcommand\Proof{\noindent{\bf Доказательство. }}
\newcommand\Endproof{\nopagebreak\strut%
\nopagebreak\hfill\nopagebreak$\Box$\medbreak}
\newcommand\Noproof{\nopagebreak\strut%
\nopagebreak\hfill\nopagebreak$\Box$\par}

\def\Fol{\operatorname{Fol}}
\def\Card{\operatorname{Card}}

\def\Ret{\operatorname{Ret}}

\def\Const{\operatorname{Const}}

\begin{document}

\title{{\small УДК 512+519.17+517.987} \\ Описание нормальных базисов
  граничных алгебр.} \author {А.~Я.~Белов, А.~Л.~Чернятьев} \date{}

\maketitle

\abstract{В работе исследуются нормальные базисы для алгебр с
 медленным ростом. Работа поддержана грантом РФФИ № 14-01-00548.}

\section {Введение.}

Для произвольной алгебры $A$ через $V_A(n)$ обозначим размерность
вектрного пространства, порожденного мономами длины не больше $n$.
Пусть $T_A(n)=V_A(n)-V_A(n-1)$. Если алгебра однородна, то $T_A(n)$
есть размерность векторного пространства, порожденного мономами
длины ровно $n$.

Известно (см.) что либо $\lim_{n\to
\infty}{(T_A(n)-n)}=-\infty$ (в этом случае есть альтернатива ({\bf
Bergman Gap Theorem}): либо $\lim{V_A(n)}=C<\infty$ и тогда $\dim A
< \infty$, либо $V_A(n)=O(n)$ и алгебра имеет {\bf медленный рост}),
либо  $T_A(n)-n < \Const$, либо, наконец, $\lim_{n\to
\infty}{(T_A(n)-n)}=\infty$.

В последнем случае для любой функции $\phi(n)\to \infty$ и любой
$\psi(n)=e^{o(n)}$ существует алгебра $A$ такая, что для
бесконечного множества натруальных чисел $n\in L\subset {\mathbb N}\
T_A(n)>\psi(n)$ и для бесконечного множества натуральных чисел $n\in
M\subset {\mathbb N}\ T_A(n)<n+\phi(n)$ (см.).

Таким образом, в случае $\lim_{n\to \infty}{(T_A(n)-n)}=\infty$ рост
может быть хаотичным и мы этот случай не рассматриваем. Случай,
когда $T_A(n)-n < \Const$ (т.е. случай алгебр {\bf медленного
роста}) исследовался Дж.~Бергманом и Л.~Смоллом. Нормальные базисы
для таких алгебр исследованы в работе \cite{BBL}. Назовем алгебру
{\bf граничной}, если $T_A(n)-n<\Const$. Наша цель состоит в
описании нормальных базисов граничных алгебр.

Прежде всего отметим, что такое описание сводится к мономиальному
случаю. В самом деле, пусть $a_1,\ldots,a_s$ -- образующие алгебры
$A$. Порядок $a_1\prec\ldots\prec a_s$ индуцирует порядок на
множестве слов алгебры $A$ (сперва по длине, затем
лексикографически). Назовем слово {\bf неуменьшаемым} (или {\bf
неприводимым}), если его нельзя представить в виде линейной
комбинации меньших слов. Множество неуменьшаемых слов образует {\bf
нормальный базис} алгебры $A$ как векторного пространства.
Рассмотрим фактор свободной алгебры $k<\hat{a}_1,\ldots,\hat{a}_s>$
($k$~-- основное поле) по множеству слов из нормального базиса.
Получится мономиальная алгебра с тем же нормальным базисом и, стало
быть, с той же функцией роста.

Назовем {\bf обструкцией} приводимое (то есть уменьшаемое) слово $u$
такое, что любое его подслово неприводимо. Сверхсловом в алгебре $A$
(правым, левым, двусторонним) называется сверхслово $W$ такое, что
любое его конечное подслово ненулевое. Аналогично определяется {\bf
неприводимое сверхслово} в алгебре $A$.

Для описания нормальных базисов граничных алгебр нам потребуется ряд
определений и утверждений, используемых в комбинаторике слов.

\section {Комбинаторика слов  и графы Рози.}

\subsection{Пространство слов.}
Здесь и далее через  $A$ будем обозначать конечный алфавит, то есть
непустое множество элементов (символов).  Через $A^+$ обозначим множество
всех конечных последовательностей, символов, или {\bf слов}.

Конечное слово всегда может быть единственным образом
представлено в виде

$w = w_1 \cdots w_n$, где $w_i \in A, 1\leq i \leq n$. Число $n$
называется {\bf длиной} слова $w$ и обозначается $|w|$

Множество $A^+$ всех конечных слов над $A$ образует простую
полугруппу, где полугрупповая операция определяется как
конкатенация (приписывание).

Если к множеству слов добавить элемент $\Lambda$ (пустое слово),
то получим свободный моноид $A^*$ над $A$. Длина $\Lambda$ по
определению равна 0.

Слово $u$ есть {\bf подслово} (или {\bf фактор}) слова $w$, если
существуют слова $p,q\in A^+$ такие, что $w = puq$.

Если слово $p$ (или $q$) равно $\Lambda$, то $u$ называется {\bf
префиксом} (или {\bf суффиксом}) слова $w$.

\begin{definition}
Слово $W$ называется {\bf рекуррентным}, если каждое его подслово
встречается в нем бесконечно много раз (в случае двустороннего
бесконечного слова, каждое подслово встречается бесконечно много
раз в обоих направлениях). Слово $W$ называется {\bf
равномерно-рекуррентным} или ({\bf р.р словом}), если оно
рекуррентно и для каждого подслова $v$ существует натуральное
$N(v)$, такое, что для любого подслова $W$ $u$ длины не менее, чем
$N(v)$, $v$ является подсловом $u$.
\end{definition}

Пусть $W$ -- бесконечное слово. Для любого подслова $v$ можно
определить множество {\bf возвращаемых} слов $\Ret_W(v)$, а
именно, слово $u$ -- {\bf возвращаемое} для $v$, если $vuv$ --
подслово $W$ и $v$ -- не является подсловом $u$. Ясно, что для
рекуррентных слов множество возвращаемых слов $\Ret(v)$ каждого
подслова $v$ будет непустым, а в случае равномерной рекуррентности
множество длин слов из $\Ret(v)$ будет ограниченным.

\begin{definition}
Для произвольного слова $W$ можно определить функцию роста:
$$
T_W(n)=\Card F_n(W)
$$
\end{definition}
Ясно, что если $T_W(n)=0$ для какого-то $n$, то $W$ -- конечное
слово. В противном случае бесконечное.

\begin{definition}
Рассмотрим бесконечное (одностороннее или двустороннее) слово над
алфавитом $A$. Пусть $v$ -- его подслово и $x\in A$. Тогда
\begin{enumerate}
\item Символ $x$ -- левое (правое) расширение $v$, если $xv$
(соотв. $vx$) принадлежит $F(W)$.

\item Подслово $v$ называется  левым (правым) специальным
подсловом, если для него существуют два или более левых (правых)
расширения.

\item Подслово $v$ называется биспециальным, если оно является и
левым, и правым специальнм подсловом одновременно.

\item Количество различных левых (правых) расширений подслова
назовем  левой (правой) валентностью этого подслова.

\end{enumerate}
\end{definition}

\subsection{Графы Рози подслов.}

Удобным инструментом для описания слова $W$ является {\bf графы
подслов}, или графы Рози (Rauzy's graphs), введенные Рози
которые строятся следующим образом: {\bf $k$-граф}
слова $W$ -- ориентированный граф, вершины которого
взаимнооднозначно соответствуют подсловам длины $k$ слова $W$, из
вершины $A$ в вершину $B$ ведет стрелка, если в $W$ есть подслово
длины $k+1$, у которого первые $k$ символов -- подслово
соответствующее $A$, последние $k$ символов -- подслово,
соответствующее $B$. таким образом, ребра $k$-графа биективно
соответствуют ($k+1$)-подсловам слова $W$.

Ясно, в $k$-графе $G$ слова $W$ правым специальным словам
соответствуют вершины, из которых выходит (соотв. в которые
входит) больше одной стрелки. Такие вершины мы будем называть
развилками. Граф $G$ будем называть {\bf сильно связным}, если из
любой вершины в любую вершину можно перейти по стрелкам.

{\bf Последователем} ориентированного графа $G$ будем называть
ориентированный граф $\Fol(G)$ построенный следующим образом:
вершины графа $G$ биективно соответствуют ребрам графа $G$, из
вершины $A$ в вершину $B$ ведет стрелка, если в графе $G$ конечная
вершина ребра $A$ является начальной вершиной ребра $B$.

Связность графов Рози и рекуррентность соответсвующего слова
связаны естественным образом. Имеет место следующее
\begin{proposition}\label{reccur}
Пусть $W$ -- бесконечное (в одну сторону) слово. Следующие условия
эквивалентны:

\begin{enumerate}
 \item Слово $W$ -- рекуррентно.

 \item Для всех $k$ соответствующий
$k$-граф слова $W$ является сильно связным.

\item Каждое подслово $W$ встречается не меньше двух раз.

\item Любое подслово является продолжаемым слева.
\end{enumerate}
\end{proposition}

В терминах графов Рози можно выразить такие важные понятия, как
рост подслов, множество запрещенных подслов, минимальные
запрещенные слова и т.д.

\section{Слова Штурма.}
Хорошо изученным примером равномерно-рекуррентных слов с медленным
ростом являются слова Штурма. Дадим несколько определений.

\begin{definition}
Бесконечное вправо (соответственно, двусторонее бесконечное) слово
$W=(w)_{n\in \mathbb{N}}$ (соответственно, $W=(w)_{n\in
\mathbb{Z}}$ ) над конечным алфавитом $A$ называется словом
Штурма, если его функция сложности равна $T_W(n)=n+1$ для всех
$n\geq 0$.
\end{definition}

{\bf Примечание 1.} Так как $T_W(1)=2$, то слово Штурма по
определению является словом над двухсимвольным алфавитом.

\begin{definition}
Слово $W$ называется сбалансированным, если для любых его двух
подслов $u,v$ одинаковой длины выполняется неравенство:

\begin{equation} \label{bal}
||u|_a-|v|_a| \leq 1
\end{equation}
\end{definition}

Рассмотрим еще один способ конструирования бесконечных сверхслов над
бинарным алфавитом.

Пусть $M$ -- компактное метрическое пространство, $U\subset M$~--- его открытое подмножество, $f:M\to M$ -- гомеоморфизм компакта в
себя и $x\in M$ -- начальная точка.

По последовательности итераций можно построить бесконечное слово
над бинарным алфавитом:

$$
w_n=\left\{
\begin{array}{rcl}
   a,\ f^{(n)}(x_0)\in U\\
   b,\ f^{(n)}(x_0)\not\in U\\
\end{array}\right.
$$
которое называется {\bf эволюцией} точки $x_0$.

Для слов над алфавитом, состоящим из большего числа символов нужно
рассмотреть несколько характеристических множеств: $U_1,\ldots
,U_n$.

Пусть
$\mathbb{S}^1$ -- окружность единичной длины, $U\subset\mathbb{S}^1$
-- дуга длины $\alpha$, $T_\alpha$ -- сдвиг окружности на
иррациональную величину $\alpha$.
Слова, получаемые такими динамическими системами, иногда называют
{\bf механическими словами} ({\bf mechanical words}).

Замечательным фактом, описываающим структуру слов Штурма является {\bf
  Теорема эквивалентности.}

\begin{theorem}\label{TheorEq}
Следующие условия на слово $W$ эквивалентны:
\begin{enumerate}
\item Слово $W$ имеет функцию сложности $T_W(n)=n+1$).

\item Слово $W$ сбалансированно и не периодично.

\item Слово $W$ является механическим, то есть порождается системой
  $(\mathbb{S}^1, U, T_\alpha)$ с иррациональным $\alpha$.
\end{enumerate}
\end{theorem}
\Noproof

\section{Слова с минимальной функцией роста.}\label{graph}

Естественным обобщением слов Штурма являются  слова с минимальной функцией роста, то есть с функцией роста,
удовлетворяющей соотношению $F_W(n+1)-F_W(n)=1$ при всех
достаточно больших $n$.

В этой части мы построим динамическую систему,
которая аналогичным образом могла бы порождать слова медленного
роста. В этой части мы будем рассматривать слова над произвольным
алфавитом $A=\{a_1,a_2,\ldots ,a_n\}$.

Пусть слово $W$ будет словом медленного роста, то есть, начиная с
некоторого $N$ верно $T_W(k+1)=T_W(k)+1, k\geq N$.

Рассмотрим $k$-графы Рози этого слова для $k \geq N$. Так как ребра
графа соответствуют ($k+1$)-словам, то в этом графе $m$ вершин и
$m+1$ ребро и, соответственно, имеет  одну входящую и одну
выходящую развилку, которые, возможно, совпадают (если данная
вершина соответствует биспециальному слову).

Имеется два типа таких графов:

\vbox{
\begin{center}
\includegraphics{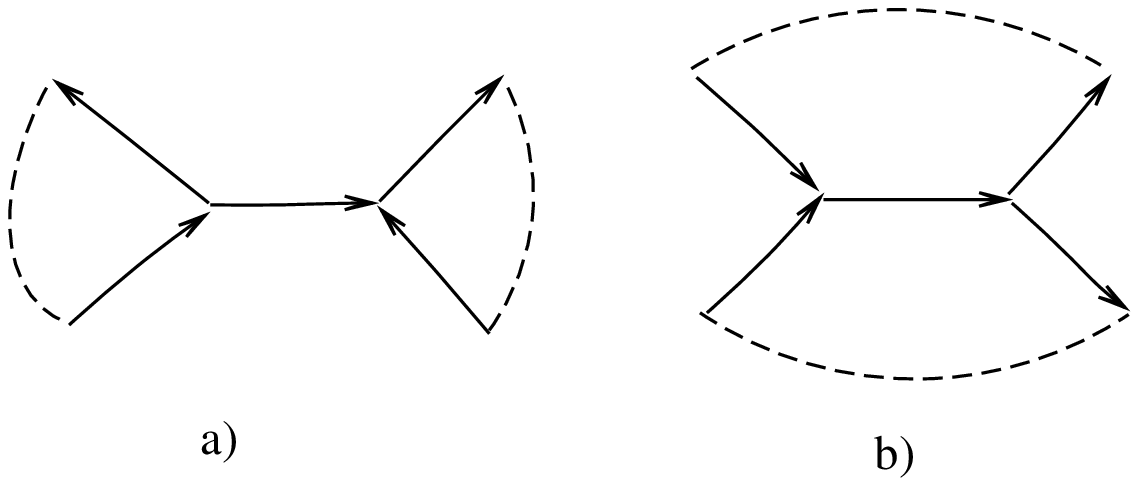}\\
{\bf Два типа графов с одной входящей и одной входящей развилкой.}
\end{center}
}\label{RisGraphRauzy1}

Легко видеть, что граф типа a) не является сильно связным, а в
случае b) -- является. Слова с $k$-графами типа b) имеют вид
$u^{\infty}Vv^{\infty}$  и не являются равномерно-рекуррентными

Нас интересует случай a).

Назовем путь, ведущий из входящей развилки в выходящую {\bf
перегородкой}, а два пути из выходящей развилки будем называть
{\bf дугами}.

\begin{proposition}
Пусть $k$-граф слова $W$ имеет перегородку длины $l\geq 1$ и дуги
длины $r,s$. Тогда последователь $\Fol(G)$ имеет перегородку длины
$l-1$, дуги длины $r+1$ и $s+1$ и совпадает с $(k+1)$-графом слова
$W$.

\end{proposition}

Рассмотрим теперь предельный случай, когда перегородка
вырождается, то есть входящая развилка совпадает с выходящей. В
этом случае развилка соответствует биспециальному слову.

\begin{proposition}
Последователь графа с вырожденной перегородкой имеет две входящие
и две выходящие развилки

\end{proposition}

Пусть $u$ -- биспециальное подслово $W$, то есть при некоторых
$a_i,a_j,a_r,a_s\in A$ $a_i u,a_j u, ua_r, ua_s$ -- тоже подслова
$W$. Тогда развилками (входящими и выходящими) в $\Fol(G)$ будут
вершины соответствующие этим словам.

Таким  образом, для $(k+1)$-графа слова $W$ имеется $4$
возможности для удаления одного ребра из $\Fol(G)$,
соответствующего минимальному запрещенному $(k+2)$-слову: $a_i u
a_r, a_i u a_s, a_j u a_s, a_j u a_r$, В двух случаях мы получаем
сильно связный граф, а в двух -- не сильно связный.

Непосредственной проверкой доказывается
\begin{proposition}
Пусть в графе $G$ перегородка вырождена, а дуги имеют длину $r,s$
соответственно. Тогда $(k+1)$-граф имеет вид $(s-1,1,r+1)$ или
$(r-1,1,s+1)$.
\end{proposition}

\begin{proposition}\label{graph:3}
Любой граф $G$ с одной входящей и одной выходящей развилкой имеет
предшественника, причем только одного.
\end{proposition}

\Proof Пусть граф имеет вид $(l,r,s)$. Заметим, что графа с
$r=s=1$ не бывает. В случае, если $r,s>1$ граф предшественника
имеет вид $(l+1,r-1,s-1)$, если имеет вид $(l,1,s),s>1$ , то
предшественник -- $(0,l+1,s-1)$, если $(l,r,1),r>1$ -- то
$(0,s+1,l-1)$. Граф вида $(0,1,k)$ имеет предшественника
$(0,1,k-1)$.\Endproof

Ясно, что любой граф может быть получен из графа типа $(0,1,2)$
или $(0,2,1)$. Из предложения \ref{graph:3} следует, что для слова
с минимальной функцией роста $W$ существует такое слово Штурма $V$
и такие натуральные $n$ и $l$, что  $k$-графы для слова $W$
совпадают с $(k+l)$-графами слова Штурма $V$ для $k\geq n$.

\begin{proposition}
Пусть $W$ -- слово с минимальной функцией роста, тогда существует
натуральное $k$, что для любого подслова $v\subset W$ такого, что
$|v|\geq k$, сущствует ровно два возвращаемых слова.
\end{proposition}

\Proof Рассмотрим слово Штурма, такое что эволюция $k$-графов Рози
слова $W$ и графов слова $V$ совпадают, начиная с некоторого $n$.
Тогда $k$-словам слова $W$ биективно соответствуют $(k+l)$-словам
слова $V$. Поскольку для подслов слова Штурма существует ровно $2$
возвращаемых слова, то это же верно и для подслов слова
$W$.\Endproof

Из этих же соображений доказывается
\begin{proposition}
рекуррентное слово с минимальной функцией роста равномерно
рекуррентно.
\end{proposition}
\Noproof

\subsection{Конструкция динамической системы}

Теперь мы покажем, как по соответствующему графу слов построить
динамическую систему, которая порождала бы данное слово.

Построение динамической системы мы осуществим двумя способами.
Первый способ опирается на уже известный результат о словах
Штурма. Второй способ -- более конструктивный и обобщающийся на
случай нескольких развилок.

Итак, пусть слово $W$ имеет минимальный рост, т.е., с
некоторого $k$ его $k$-графы слова $W$ имеют вид a).

Поскольку каждый граф имеет единственного предшественника, то
существует последовательность графов типа a) $G'_1, G'_2, \ldots,
G'_n$ таких, что $G'_{l+1}$ является предшественником $G'_l$, а
граф $G'_n$ совпадает с $k$-графом слова $W$.

Иными словами, существует эволюция $k$-графов типа a), начиная с
$k=1$, что $n$-ый граф в эволюции совпадает с $k$-графом слова
$W$, ($n+1$)-ый  совпадает с ($k+1$)-графом слова $W$ и т.д.

Последовательность графов такого типа однозначно соответствует
некоторому слову Штурма $V$.

Значит, существует взаимнооднозначное соответствие между
$k$-словами слова $W$  и $n$-словами слова Штурма $V$, при этом
соответствие продолжается на ($k+1$)-слова $W$ и ($n+1$)-слова $V$
и т.д.

Разберем сначала случай, когда $k=1$, то есть уже $1$-граф слова
$W$ имеет вид а).  Символам алфавита $A$ взаимнооднозначно
соответствуют $n$-слова некоторого слова Штурма $V$. По теореме
эквивалентности, слово $V$ порождается сдвигом окружности
$T_\alpha$. При этом $n$-словам в
динамике соответствуют интервалы разбиения: слову $w=w_1w_2\ldots
w_n$ ($w_i\in \{a,b\}$) соответствует интервал
$$
I_w=T^{(-n+1)}_\alpha(I_{w_1})\cup T^{(-n+2)}_\alpha(I_{w_2}) \cup
\ldots \cup I_{w_n},$$ где $I_{w_i}$ -- характеристичесий интервал
для $w_i$.

Для сдвига окружности интервалы $I_w$, соответствующие $n$-словам,
будут иметь вид $I_w=(n_i\alpha,n_{i+1}\alpha)$, где $n_i$ --
некоторые целые числа.

Тогда характеристическим множеством для каждого символа из
алфавита $A$ будет являться интервал $I_w$, такой, что слово $w$
соответствует данному символу.  Ясно тогда, что слово $W$ будет
порождаться тем же сдвигом $T_\alpha$ и характеристическими
множествами $I_w$.

Теперь разберем общий случай. Пусть $k$-словам слова $W$
соответствуют $n$-слова слова Штурма $V$.  Точно также мы можем
построить соответствие между интервалами разбиения для $n$-слов
слова Штурма $I_w$ и $k$-словами слова $V$.  Построим
характеристическое множество для каждого символа из $A$.

А, именно, символу $a_i$ поставим в соответствие все интервалы,
которые соответствуют словам, начинающимся на $a_i$.  В этом случае
характеристическое множество для произвольного символа может быть
несвязным и представлять собой объединение нескольких интервалов.

Покажем, что при таком выборе характеристических множеств
найдется точка, чья эволюция будет совпадать с $W$.

Действительно, пусть слово $W$ начинается с некоторого $k$-слова
$w=w_1w_2\cdots w_k$. Ему мы поставим в соответствие интервал $I_w$.
Рассмотрим следующий $k+1$ символ $w_{k+1}$; этому ($k+1$)-слову мы
поставим в соответствие интервал $I_w\cup I_{w'}$, где $w'=w_2w_3
\cdots w_{k+1}$, и т.д.

Таким образом, доказана следующая

\begin{theorem}\label{MinGrow}
Пусть $W$ -- рекуррентное слово над произвольным конечным
алфавитом $A$. Тогда следующие условия на слово $W$ эквивалентны:
\begin{enumerate}

\item Существует такое натуральное $N$, что функция сложности
слова $W$ равна $T_W(n)=n+K$, для $n\geq N$  и некоторого
постоянного натурального $K$.

\item Существуют такое иррациональное $\alpha$ и целые $n_1,n_2,
\ldots,n_m$, что слово $W$ порождается динамической системой
$$
(\mathbb{S}^1,T_\alpha, I_{a_1},I_{a_2},\ldots ,I_{a_n},x),
$$
где
$T\alpha$ -- сдвиг окружности на иррациональную величину $\alpha$,
$I_{a_i}$ -- объединение дуг вида $(n_j\alpha,n_{j+1}\alpha)$.

\end{enumerate}

\end{theorem}

\section{Основная теорема.}

\begin{theorem}[Описание нормальных базисов граничных алгебр]\label{ThMainLast}
Пусть $A$ -- граничная алгебра, $a_1,\ldots, a_s$ -- ее образующие.
Тогда имеют место два случая, каждый из которых имеет свое описание.

{\bf Случай 1.} Алгебра $A$ не содержит равномерно-рекуррентного
непериодического сверхслова. В этом случае нормальный базис алгебры
$A$ состоит из множества подслов следующего множества слов:

\begin{enumerate}
\item Одно слово вида $W=u^{\infty/2}cv^{\infty/2} \neq u^{\infty}$

\item Произвольный конечный набор $\mu$ конечных слов

\item Множество слов вида $u_i^{\infty/2}c_i$,\ $i=1,\ldots, r_1$

\item Множество слов вида $d_iv_i^{\infty/2}$,\ $i=1,\ldots, r_2$

\item Множество слов вида $e_j{(R_j)}^{k} f_j$,\ $k\in \mathbb{K}_j \subseteq
\mathbb{N}$, $j=1,\ldots,r_3$

\item Множество слов вида $W_{\alpha} = E_{\alpha}u^{n_\alpha}cv^{m_\alpha}
F_\alpha$. При этом существуют такое $c>0$, что для любого $k$
количество слов $W_\alpha$ длины $k$ меньше $c$.
\end{enumerate}

{\bf Случай 2.} Алгебра $A$ содержит равномерно-рекуррентное
непериодическое сверхслово $W$. В этом случае нормальный базис
алгебры $A$ состоит из множества подслов следующего семейства слов,
включающее в себя:

\begin{enumerate}
\item Некоторое равномерно рекуррентное слово $W$ с функцией роста
$T_W(n)=n+const$ для всех достаточно больших $n$. Описание таких
слов дано в теореме \ref{MinGrow}.

\item Произвольный конечный набор $\mu$ конечных слов

\item Множество слов вида $u_i^{\infty/2}c_i$, $i=1,\ldots,s_1$

\item Множество слов вида $d_iv_i^{\infty/2}$, $i=1,\ldots,s_2$

\item Множество слов вида $e_j(R_j)^{k}f_j$, $k\in \mathbb{K}_j \subseteq
\mathbb{N}$,\ $j=1,\ldots,s_3$

\item Множество слов вида $L_iO_iW_i$,\ $i=1,\ldots,k_1$. При этом
$W_i$ -- сверхслово, эквивалентное $W$ и $O_iW_i$ имеют вхождение
только одной обструкции (а, именно, $O_i$).

\item Множество слов вида $W'_jO'_jL'_j$, $j=1,\ldots,k_2$. При этом
$W'_j$ -- сверхслово, эквивалентное $W$ и $W'_jO'_j$ имеют вхождение
только одной обструкции (а, именно, $O'_j$).

\item Конечное множество серий вида: ${h^1}_iT_i{h^2}_i$,
$i=1,\ldots,s$. При этом:

a) слово $T_i$ содержит вхождение ровно двух обструкций
$O_{i_1},O_{i_2}$ (возможно, перекрывающихся)

b) Для некоторого $c>0$ $|h^1_i|+|h^2_i|<c$ при всех $i$.

с) Существует $m>0$ такое, что для любого $k$ имеется не более $k$
подслов длины $m$, вида $(8)$ и не являющихся подсловами слова $W$.

\end{enumerate}
\end{theorem}

\Proof Алгебру $A$ считаем мономиальной. Назовем слово $v$ алгебры
$A$ {\bf хорошим}, если для любого $n$ существуют сколь угодно
длинные слова $w_1,w_2$, $|w_1|>n,|w_2|>n$ такие, что $w_1vw_2$ есть
подслово алгебры $A$. Обозначим $T_{RL}(n)$  количество хороших слов
длины $n$. Известно, что если $T_{RL}(n)=T_{RL}(n+1)$ при некотором
$n$, то алгебра $A$ имеет медленный рост (см. \cite{BBL}).  В силу
граничности алгебры $A$ $T_{RL}(n)\geq T_{RL}(n+1)+1$,  при всех
достаточно больших $n$ неравенство превращается в равенство (иначе
$\lim_{n\to \infty} (T_{RL}(n)-n)= \infty$ и, как следствие
$\lim_{n\to \infty} (T_{A}(n)-n)= \infty$).

При этом граф Рози имеет развилку и, как следствие, два цикла,
эволюция графа Рози устроена следующим образом: либо граф теряет
сильную связность и имеет вид a). В этом случае его дальнейшая
эволюция однозначна, она отвечает слову вида $u^{\infty/2} c
v^{\infty/2}$ и мы имеем случай $1$, либо граф Рози все время
остается сильно связным. Тогда эволюция связной компоненты, в
которой есть развилка, асимптотически эквивалентна эволюции графа
Рози некоторого слова Штурма. Если оно имеет вид  $u^{\infty/2} c
u^{\infty/2}$, то имеет место случай $1$, иначе оно
равномерно-рекуррентно и имеет место случай $2$.

В случае 1 все обструкции для слова $u^{\infty/2} c u^{\infty/2}$
имеют ограниченную длину (см. \cite{BBL}). Можно сделать следующее

\medskip
{\bf Наблюдение.}\ {\bf Количество ненулевых слов длины $n$, не
являющихся подсловами слова $u^{\infty/2}cu^{\infty/2}$ не
превосходит константы (не зависящей от $n$).}
\medskip

Напомним предложение из работы \cite{BBL}:

\begin{proposition}
Пусть $SW=WT$. Тогда $W$ имеет вид: $s^k s_1$, где $s_1$ -- начало
слова $s$.
\end{proposition}

Из данного предложения и только что сделанного наблюдения вытекает

\begin{corollary}
Пусть $|u|=k$, $|v_1|=|v_2|=l$. Тогда либо количество подслов длины
$k+m$, (где $m\le k$) слова $v_1uv_2$ не менее $m+1$, либо
$u=s^ks'$, для некоторого слова $s$, при этом $s'$ -- начало $s$ и
$v_1=v'_1s'$.
\end{corollary}

Из данного следствия и наблюдения получается

\begin{proposition}
Существует константа $K$, зависящая только от граничной алгебры $A$,
такая, что для любой обструкции $O$ в слове $W$ либо при некотором
$m$ для любых $v_1,v_2$, $|v_1|\geq m$, $|v_2|\geq m$ $v_1Ov_2$
является нулевым словом алгебры (число $m$ не зависит от выбора
обструкции), либо $|v|\leq K$.
\end{proposition}

Из данного предложения следует, что словами алгебры $A$ являются
либо слова, содержащие не более двух обструкций, причем каждая
обструкция находится на ограниченном расстоянии от одного из концов,
либо подслова слов вида $R_iu^k_iT_i$. А все такие типы слов описаны
в условии теоремы \ref{ThMainLast}. \Endproof

\end{document}